\documentclass{amsart}
\usepackage{amssymb}
\usepackage{enumerate}
\def\CC{\mathbb C}
\def\CO{\mathcal O}
\def\conv{\operatorname{conv}}
\def\phi{\varphi}
\def\Ree{\operatorname{Re}}
\def\RR{\mathbb R}
\def\eps{\varepsilon}
\def\too{\longrightarrow}
\def\Cal{\mathcal}
\def\wdht{\widehat}
\def\wdtl{\widetilde}

\makeatletter
\def\th@mytheorem{%
  \let\thm@indent\noindent
  \thm@headfont{\bfseries}
    \itshape
}
\def\th@myremark{%
  \let\thm@indent\noindent
  \thm@headfont{\bfseries}
}
\makeatother
\theoremstyle{mytheorem}
\newtheorem{Theorem}{Theorem}
\theoremstyle{myremark}
\newtheorem{Remark}[Theorem]{Remark}
\begin{document}
\title{The envelope of holomorphy of a classical truncated tube domain}

\author[M.~Jarnicki]{Marek Jarnicki}
\address{Jagiellonian University, Faculty of Mathematics and Computer Science, Institute of Mathematics,
{\L}ojasiewicza 6, 30-348 Krak\'ow, Poland}
\email{Marek.Jarnicki@im.uj.edu.pl}

\author[P.~Pflug]{Peter Pflug}
\address{Carl von Ossietzky Universit\"at Oldenburg, Institut f\"ur Mathematik,
Postfach 2503, D-26111 Oldenburg, Germany}
\email{Peter.Pflug@uni-oldenburg.de}

\begin{abstract}
We present the envelope of holomorphy of a classical truncated tube domain.
\end{abstract}



\maketitle

S.~Bochner showed that the envelope of holomorphy of a tube domain $G=D+i\RR^n\subset\CC^n$, $n\geq 2$, is given by the convex hull of $G$. Discussions of so called generalized tube domains $G=D_1+iD_2$ (also called truncated tube domains)
have shown that Bochner's result becomes false, i.e. the envelope of holomorphy of $G$ is, in general, a strict subdomain of $\conv(D_1)+i\conv(D_2)$ (see \cite{JarPfl2020}). The simplest known example of this type is given by the following domain

$$ G:=\{x+iy\in\CC^n: R_1'<\|x\|<R_1'',\;\|y\|<R_2\},\;n\geq 2,$$
where $\|\;\|$ stands for the Euclidean norm in $\RR^n$ and $0\leq R_1'<R_1''$, $0<R_2$, see \cite{JarPfl2020}.

Recently, discussing the tube theorem, (see \cite{Nog2020}) J.~Noguchi asked for the envelope of holomorphy $\Cal E(G)$ of $G$. Here we present the answer.

\begin{Theorem} $\Cal E(G)= \{x+iy\in\wdht G: \|y\|^2<\|x\|^2-(R_1'^2-R_2^2)\} $, 
where $\wdht G$ is the convex hull of $G$.  In particular, the envelope is univalent.
\end{Theorem}

\begin{Remark} This result shows that the envelope of our truncated tube domain is never the convex hull of $G$ except when $R_1'=0$ or $R_2=\infty$ (see the paper of J.~Kajiwara \cite{Kaj1963}); the last case is just  Bochner's tube theorem  while the other case can be also handled using \cite{Shi1968} if $n>2$ and \cite{Pfl1980}, \cite{Kar1991} if $n=2$. Moreover, it is univalent; but we don't know whether the envelope of every generalized tube domain is univalent.
\end{Remark}

\begin{proof}
Without loss of generality we may assume that $R_1''=1$;  then we set $R_1:=R_1'$ to simplify the notation, i.e. the domain we discuss is
$$
G:=\{x+iy\in\CC^n:R_1<\|x\|<1,\;\|y\|<R_2\}.
$$

The proof is done via induction. Let $n=2$. Then it is based on the work of S.M.~Ivashkovich (\cite {Iva1982}),  in particular on Lemma 8 in this work.

Let $f\in\CO(G)$.

Note that $G\subset \{z\in\CC^2:\Ree(z_1^2+z_2^2)-(R_1^2-R_2^2)>0\}$. For $\eps\in(0,\min\{R_2^2, 1-R_1^2\})$, $\eps\neq R_2^2-R_1^2$, put
\begin{align*}
\phi_\eps(z):&= z_1^2+z_2^2-(R_1^2-R_2^2)-\eps,\; z=(z_1,z_2)\in\CC^2, \\
\Pi_\eps:&=\{z\in\CC^2:\Ree\phi_\eps(z)=0\}, \text{and finally } \\
\Pi_\eps^+:&=\{z\in\CC^2:\Ree\phi_\eps(z)>0\}.
\end{align*}

Put $R_{1,\eps}:=\sqrt{R_1^2+\eps}$, and fix $R_{1,\eps}^\ast\in (R_1,R_{1,\eps})$.

Finally, define $D_\eps:=\{z\in\CC^2: \|x\|<R_{1,\eps}^\ast, \|y\|^2<\|x\|^2-(R_1^2-R_2^2)-\eps\}\subset \Pi_{\eps}^+$.

Then $S_\eps:=\overline{\partial D_\eps\setminus(\overline {D_\eps}\cap\Pi_\eps)}$ is connected and $S_\eps\subset G$. Note that $f$ is holomorphic in a neighborhood of $S_\eps$. Then, using Lemma 8 of \cite{Iva1982},  $f$ extends to a holomorphic function $f_\eps$ on $D_\eps$ which is identically equal to $f$ near $S_\eps$. Using the identity theorem and $\eps\too 0$ we end up with a holomorphic function on $\Cal E(G)$.

Observe that $\Cal E(G)$  is the intersection of a convex domain and a Levi-flat pseudoconvex domain and so it is pseudoconvex. Therefore we have got the envelope of holomorphy of $G$.

Now let $n\geq 3$ and assume that the theorem is true in dimension $n-1$. Take an $f\in\CO(G)$ and fix a point $z_n=x_n+i y_n$ with $|x_n|<1$, $|y_n|<R_2$. Then $G$ intersects the hyperplane $\CC^{n-1}\times\{z_n\}$
along
$$
G(z_n):= \{\wdtl z\in\CC^{n-1}:R_1^2-|x_n|^2<\|\wdtl x\|^2<1-|x_n|^2, \;\|\wdtl y\|^2<R_2^2-|y_n|^2\}.
$$
If $|x_n|>R_1$, then $G(z_n)$ is convex and therefore, pseudoconvex; in particular, $\Cal E(G(z_n))=G(z_n)$. On $G(z_n)$ put $\wdtl f(\cdot,z_n):=f(\cdot,z_n)$. 

If $|x_n|\leq R_1$, then $G(z_n)$ is our domain, now in $\CC^{n-1}$.
Therefore, by induction, there exists an  $\wdtl f(\cdot, z_n)\in\CO(\Cal E(G(z_n))$ extending $f(\cdot, z_n)$. 


Summarizing, we end up with a function $\wdtl f$ on $\Cal E(G)=\bigcup_{|x_n|<1, |y_n|<R_2}\Cal E(G(z_n))$.
Recall that for all $z_n=x_n+iy_n$ with $|x_n|<1$ and $|y_n|<R_2$, the function  $\wdtl f(\cdot, z_n)$ is holomorphic. Moreover, $\wdtl f=f$ on $G$, i.e.~$\wdtl f$ is holomorphic on $G$.

It remains to show that $\wdtl f$ is holomorphic on $\Cal E(G)$. Fix a point $a+ib\in\Cal E(G)$. If $\|a\|> R_1$, then $a+ib\in G$ and so $\wdtl f$ is holomorphic in a neighborhood of $a+ib$.

Now assume that $\|a\|\leq R_1$. Observe that it suffices to find open connected sets $\varnothing\neq U'_1\subset U_1\subset\CC^{n-1}$, $U_2\subset\CC$ with $a+ib\in U_1\times U_2\subset \Cal E(G)$ and $U'_1\times U_2\subset G$. Then Hartogs Theorem on separate analyticity (see Theorem 1.1.10(a) from  \cite{JarPfl2011}) will show that $\wdtl f$ is holomorphic on $U_1\times U_2$, and therefore, in a neighborhood of $a+ib$.

 Take a positive real number $\eps$ such that $R_1^2+2\eps <1$ and
 $\|b\|^2+2\eps<\|a\|^2-(R_1^2-R_2^2)$. Put
\begin{align*} U_1:&=\{\wdtl x+i\wdtl y \in\CC^{n-1}:\|\wdtl a\|^2-\eps/4<\|\wdtl x\|^2<1-|a_n|^2-\eps/2, \|\wdtl y\|^2<\|\wdtl b\|^2+\eps/2\}, \\
U_2:&=\{x_n+iy_n\in\CC: |a_n|^2-\eps/4<|x_n|^2<|a_n|^2+\eps/2, |y_n|^2<|b_n|^2+\eps/2\},\\
U_1':&=\{\wdtl z=\wdtl x+i\wdtl y\in U_1: 1-|a_n|^2-3\eps/2<\|\wdtl x\|^2<1-|a_n|^2-\eps/2\}.
\end{align*}
Then $a+ib\in U_1\times U_2\subset\Cal E(G)$ and $\varnothing\neq U_1'\times U_2\subset G$. Therefore, using Hartogs theorem, $\wdtl f$ is holomorphic in a neighborhood of $a+ib$.

Hence $\wdtl f$ is holomorphic on $\Cal E(G)$ and we have shown that $\Cal E(G)$ is a holomorphic extension of $G$. On the other hand, as above, it is easy to see that $\Cal E(G)$ is pseudoconvex. Therefore, it is the envelope of holomorphy of $G$.
\end{proof}

\begin{Remark}  Observe that the generalized tube domain $G=D_1+iD_2$ discussed above has a special symmetry; namely, whenever $x,x'\in\RR^n$, $\|x\|=\|x'\|$, and $x\in D_1$, then $x'\in D_1$ (resp. whenever $y,y'\in\RR^n$, $\|y\|=\|y'\|$, and $y\in D_2$, then $y'\in D_2)$. Moreover $D_2$ is convex. So one may ask whether any generalized tube domain with this kind of an additional geometric property (e.g. substitute in the definition of $G$ the norm $\|\cdot\|$ by another $\RR$-norm) has a univalent envelope of holomorphy and what is its description.
\end{Remark}

Let us end this short note with emphasizing some  open general problems (see \cite{Nog2020}):
\begin{Remark}

a) Does every generalized tube domain have a univalent envelope of holomorphy ?

b) In case of univalence what is the description of the envelope ?

c) What is the shape of the envelope of holomorphy in general ?
\end{Remark}

Acknowledgments.  The authors thank Professor H.~Boas, Professor J.~Noguchi, and the referee
for many helpful remarks that essentially improved the presentation of this note.

\bibliographystyle{amsplain}

\end{document}